\documentclass[12pt]{amsart}
\usepackage[colorlinks=true,citecolor=black,linkcolor=black,urlcolor=blue]{hyperref}
\usepackage{amsmath}
\usepackage[english,  activeacute]{babel}
\usepackage[utf8]{inputenc}
\usepackage{amssymb}
\usepackage{amsthm}
\usepackage{graphics,graphicx}
\usepackage{enumerate}
\usepackage{array}
\usepackage{bm}
\usepackage{a4wide}
\usepackage{color, url}
\usepackage{float}

\newcommand{\seqnum}[1]{\href{http://oeis.org/#1}{\underline{#1}}}

\theoremstyle{plain}
\newtheorem{theorem}{Theorem}[section]

\newtheorem{coro}[theorem]{Corollary}

\theoremstyle{definition}

\newcommand{\Ar}{{\mathcal{A}}}
\newcommand{\E}{{\mathbb E}}
\def\last{{\textsf{last}}}
\def\parts{{\textsf{parts}}}
\newcommand{\Arr}{{\sf{Ar}}}

\title{Arndt compositions: a generating functions approach}

\date{\today}

\begin{document}

\author[D. F. Checa]{Daniel  F. Checa}
\address{Departamento de Matem\'aticas,  Universidad Nacional de Colombia,  Bogot\'a, Colombia}
\email{dcheca@unal.edu.co}
\author[J. L. Ram\'{\i}rez]{Jos\'e L. Ram\'{\i}rez}
\address{Departamento de Matem\'aticas,  Universidad Nacional de Colombia,  Bogot\'a, Colombia}
\email{jlramirezr@unal.edu.co}

\begin{abstract}
We use generating functions to enumerate Arndt compositions, that is, integer compositions where there is a descent between every second pair of parts, starting with the first and second part, and so on. In 2013, Jörg Arndt noted that this family of compositions is counted by the Fibonacci sequence. We provide an approach that is purely based on generating functions to prove this observation. We also enumerate these compositions with respect to the number of parts and the last part. From this approach, we can generalize some recent results given by Hopkins and Tangboonduangjit in 2023. Finally, we study some possible generalizations of this counting problem. 
\end{abstract}

\maketitle

\section{Introduction and Notation}

A \emph{composition} of a positive integer $n$ is a sequence of positive integers \linebreak $\sigma=(\sigma_1, \sigma_2,\ldots, \sigma_\ell)$ such that $\sigma_1+\sigma_2+\cdots +\sigma_\ell=n$. The summands $\sigma_i$ are called \emph{parts} of the composition, and $n$ is referred to as the \emph{weight} of $\sigma$ and is denoted by $|\sigma|$.  For example, the compositions of $4$ are
$$(\texttt{4}), \quad (\texttt{3},\texttt{1}), \quad(\texttt{1},\texttt{3}), \quad(\texttt{2},\texttt{2}), \quad (\texttt{2},\texttt{1},\texttt{1}), \quad (\texttt{1},\texttt{2},\texttt{1}), \quad (\texttt{1},\texttt{1},\texttt{2}), \quad (\texttt{1},\texttt{1},\texttt{1},\texttt{1}).$$

Recently, Hopkins and Tangboonduangjit \cite{HT1, HT2} have begun a combinatorial study of an interesting family of compositions with a restriction over pairwise descending parts. Specifically, an \emph{Arndt composition}  of  $n$ is a composition $\sigma=(\sigma_1, \sigma_2,\ldots, \sigma_\ell)$ of weight $n$ such that $\sigma_{2i-1}>\sigma_{2i}$ for each positive integer $i$. If the number of parts is odd, then the last inequality is vacuously true.  Let $\Ar(n)$ denote the set of Arndt compositions of weight $n$ and $\Ar=\bigcup_{n\geq 0}\Ar(n)$. For example, the Arndt compositions of $6$ are
$$(\texttt{6}), \quad (\texttt{5},\texttt{1}), \quad (\texttt{4},\texttt{2}), \quad  (\texttt{4},\texttt{1},\texttt{1}), \quad  (\texttt{3},\texttt{2},\texttt{1}), \quad  (\texttt{3},\texttt{1},\texttt{2}), \quad  (\texttt{2},\texttt{1},\texttt{3}), \quad  (\texttt{2},\texttt{1},\texttt{2},\texttt{1}).$$

This name is in honor of Jörg Arndt, who observed that the cardinality of $\Ar(n)$  is given by the $n$th Fibonacci number $F_n$ \cite{HT1}. Recall that Fibonacci numbers are defined by the recurrence relation $F_n=F_{n-1} + F_{n-2}$ for $n\geq 2$, with the initial conditions $F_0=0$ and $F_1=1$.

The goal of this paper is to enumerate Arndt compositions and some generalizations by means of generating functions.  We use (ordinary) generating functions to obtain explicit combinatorial formulas for the counting sequences. We generalize many of the results given by Hopkins and Tangboonduangjit in \cite{HT1, HT2}.  Additionally, we provide an asymptotic approximation for the expected number of parts and the last part. We also establish a connection between the number of Arndt compositions and the reduced anti-palindromic compositions with a fixed number of parts. This last family of compositions was introduced in 2022 by Andrews, Just, and Simay \cite{AJS}.  In the last two sections, we study two generalizations of Arndt compositions.

\section{Enumeration of Arndt compositions}

\begin{figure}[h]
\centering
  \includegraphics[scale=0.8]{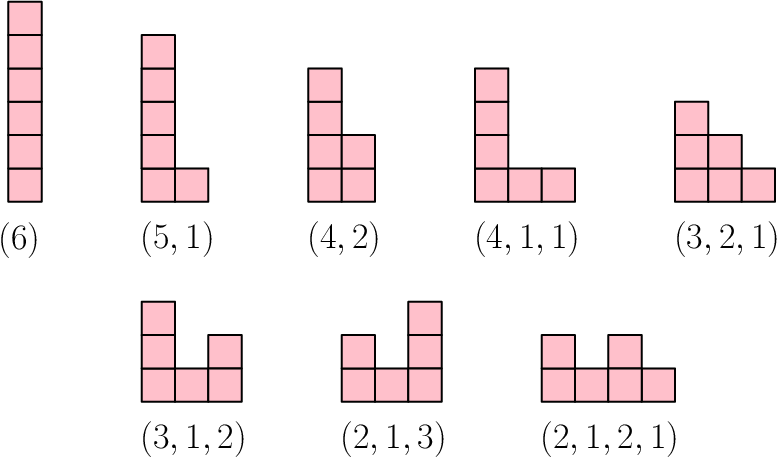}
  \caption{Arndt compositions of $n=6$.}
  \label{Fig1}
\end{figure}

A composition $(\sigma_1, \sigma_2,\ldots, \sigma_\ell)$ of weight $n$ can be represented as a \emph{bargraph} of $\ell$ columns, such that the $i$-th column contains $\sigma_i$ cells for $1\leq i \leq \ell$.  For example, in Figure \ref{Fig1} we show the Arndt compositions of $n=6$ with their bargraph representations.

We denote by $\parts(\sigma)$ the number of parts in a composition $\sigma$. We now introduce a bivariate generating function to count the number of Arndt compositions with respect to the weight and number of parts:
$$A(x,y):=\sum_{\sigma\in \Ar}x^{|\sigma|}y^{\parts(\sigma)}.$$
Note that the coefficient $x^ny^m$ of $A(x,y)$ is equal to the number of Arndt compositions of $n$ with $m$ parts. In Theorem \ref{gfun1} we give a rational generating function for $A(x,y)$.
\begin{theorem}\label{gfun1}
The generating function for Arndt compositions  with respect to the number of
parts and weight is given by
$$A(x,y)=\frac{1 - x - x^2 + x^3 + x y - x^3 y}{1 - x - x^2 + x^3 - x^3 y^2}.$$
\end{theorem}
\begin{proof}
Let $\sigma=(\sigma_1, \sigma_2)$  be an Arndt composition with two parts. From the definition we have the condition $\sigma_1>\sigma_2\geq 1$, see Figure \ref{deco1}  for a graphical representation of this case.
\begin{figure}[H]
\centering
  \includegraphics[scale=0.6]{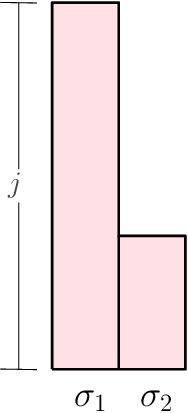}
  \caption{Decomposition of an Arndt composition with two parts.}
  \label{deco1}
\end{figure}
If $\sigma_1=j\geq 2$, then  the bivariate generating function for this case is  
$$\sum_{j\geq 2}x^{j}y\left(x+x^2+\cdots +x^{j-1}\right)y=y^2\sum_{j\geq 2}x^j\frac{x-x^j}{1-x}=\frac{x^3y^2}{(1-x)^2(1+x)}.$$
Arndt compositions (bargraphs) are the concatenation of pairs of parts (columns) as in Figure \ref{deco1}. Therefore the generating function for Arndt compositions with an even number of parts (see Figure \ref{deco1b}) is given by
 $$\sum_{m\geq 0}\left(\frac{x^3}{(1-x)^2(1+x)}\right)^my^{2m}=\frac{(1 - x)^2 (1 + x)}{1 - x - x^2 + x^3 - x^3 y^2}.$$
 \begin{figure}[H]
\centering
  \includegraphics[scale=0.6]{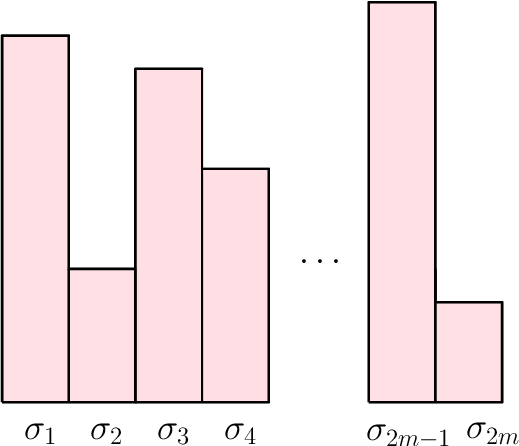}
  \caption{Decomposition of an Arndt composition with an even number of parts.}
  \label{deco1b}
\end{figure}
Analogously, if the number of parts is odd, then generating functions is given by
 $$\sum_{m\geq 0}\left(\frac{x^3}{(1-x)^2(1+x)}\right)^my^{2m}\frac{xy}{1-x}=\frac{x (1 - x^2) y}{1 - x - x^2 + x^3 - x^3 y^2}.$$
 Notice that the generating function $xy/(1-x)$ corresponds to the last part (last column). Adding the last two equations, we obtain the desired result.
\end{proof}

As a series expansion, the generating function $A(x,y)$ begins with
\begin{multline*}
A(x,y)=1+x y+x^2 y+x^3 \left(y^2+y\right)+x^4
   \left(y^3+y^2+y\right)\\+x^5 \left(2
   y^3+2 y^2+y\right)+\bm{x^6 \left(y^4+4
   y^3+2 y^2+y\right)}+O\left(x^7\right).
   \end{multline*}
Figure \ref{Fig1} shows the Arndt compositions corresponding to the bold coefficient in the above series.

Let $a(n)$ and $a(n,m)$ denote the number of Arndt compositions of $n$ and the number of Arndt compositions of $n$ with exactly $m$ parts, respectively. It is clear that $a(n)=\sum_{m\geq 1}a(n,m)$. In Table \ref{tab1} we show the first values of the sequence $a(n,m)$. This array corresponds to the entry \seqnum{A354787} in the OIES \cite{OEIS}. This array also counts the number of reduced anti-palindromic compositions, we prove this relation in Theorem \ref{bijection}.

\begin{table}[htp]
\begin{center}
\begin{tabular}{c|ccccccccccccc}
$n\backslash m$ &0& 1 & 2 & 3 & 4 & 5 & 6 & 7     \\ \hline
 0&1 &  &  &  &  &  &  &  \\
1& 0 & 1 &  &  & &  &  &  \\
2& 0 & 1 &  &  &  &  &  &  \\
3& 0 & 1 & 1 &  &  &  &  &  \\
4& 0 & 1 & 1 & 1 &  &  &  &  \\
5& 0 & 1 & 2 & 2 &  & &  &  \\
6& 0 & 1 & 2 & 4 & 1 &  &  &  \\
7& 0 & 1 & 3 & 6 & 2 & 1 &  &  \\
8& 0 & 1 & 3 & 9 & 5 & 3 &  &  \\
9& 0 & 1 & 4 & 12 & 8 & 8 & 1 &  \\
10& 0 & 1 & 4 & 16 & 14 & 16 & 3 & 1 \\
\end{tabular}
\end{center}
\caption{Values of $a(n,m)$, for $0\leq n \leq 10$ and $1\leq m \leq 7$.}
\label{tab1}
\end{table}

In Theorems \ref{teo2} and \ref{teo2b} we give combinatorial sums and recurrence relations to calculate the sequence $a(n,m)$.
\begin{theorem}\label{teo2}
For all $n, m\geq 0$ we have 
$$a(n,m)=\sum_{\ell=0}^{n-m-\lfloor\frac{m}{2}\rfloor}\binom{m + \ell - 1}{\ell}\binom{n - m - \ell - 1}{n - m - \lfloor\frac{m}{2}\rfloor - \ell}(-1)^{n - m  -\lfloor\frac{m}{2}\rfloor - \ell}.$$
Moreover, for all $n\geq 3$ and $m\geq 2$, 
$$a(n,m)=a(n-1,m)+a(n-2,m) - a(n-3,m) + a(n-3,m-2).$$
\end{theorem}
\begin{proof}
From the proof of Theorem \ref{gfun1} and the well-known identity (cf. \cite{CM})
$$\frac{1}{(1-x)^n}=\sum_{\ell\geq 0}\binom{n+\ell-1}{\ell}x^n$$ 
we have
\begin{align*}
a(n,2m)&=[x^n]\frac{x^{3m}}{(1-x)^{2m}(1+x)^m}\\
&=[x^{n-3m}]\sum_{i\geq 0}\binom{2m+i-1}{i}x^{i}\sum_{i\geq 0}\binom{m+i-1}{i}(-1)^ix^{i}\\
&=[x^{n-3m}]\sum_{i\geq 0}\sum_{\ell=0}^i\binom{2m+\ell-1}{\ell}\binom{m+i-\ell-1}{i-\ell}(-1)^{i-\ell}x^i\\
&=\sum_{\ell=0}^{n-3m}\binom{2 m + \ell - 1}{\ell}\binom{n - 2 m - \ell - 1}{n - 3 m - \ell}(-1)^{n - m - \ell}.
\end{align*}
A combinatorial formula for $a(n,2m+1)$ is obtained in a similar manner:
$$a(n,2m+1)=\sum_{\ell=0}^{n-3m-1}\binom{2 m + \ell}{\ell}\binom{n - 2 m - \ell -2}{n - 3 m - \ell-1}(-1)^{n - m - \ell-1}.$$
From these two identities we obtain the desired result. Finally, the recurrence relation follows from the equality
$(1 - x - x^2 + x^3 - x^3 y^2)A(x,y)=1 - x - x^2 + x^3 + x y - x^3 y.$
\end{proof}

\begin{theorem}\label{teo2b}
For all $n, m\geq 0$ we have 
$$a(n,m)=\sum_{\ell=0}^{\lfloor (n-m-\lfloor\frac{m}{2}\rfloor)/2\rfloor}\binom{\lfloor\frac{m}{2}\rfloor + \ell - 1}{\ell}\binom{n - 2 \lfloor\frac{m}{2}\rfloor -2\ell-1}{\lfloor\frac{m-1}{2}\rfloor}.$$
\end{theorem}
\begin{proof}
We argue similarly as in the proof of Theorem \ref{teo2}. 
\begin{align*}
a(n,2m)&=[x^n]\frac{x^{3m}}{(1-x)^{2m}(1+x)^m}=[x^{n-3m}]\frac{1}{(1-x^2)^{m}(1-x)^m}\\
&=[x^{n-3m}]\sum_{i\geq 0}\binom{m+i-1}{i}x^{2i}\sum_{i\geq 0}\binom{m+i-1}{i}x^{i}\\
&=[x^{n-3m}]\sum_{\ell \geq 0}\sum_{i\geq 0}\binom{m+\ell-1}{\ell}\binom{m+i-1}{i}x^{i+2\ell}.
\end{align*}
This with $t = i+2\ell$ implies
\begin{align*}
a(n,2m)&=[x^{n-3m}]\sum_{\ell \geq 0}\sum_{t\geq 2\ell}\binom{m+\ell-1}{\ell}\binom{m+t-2\ell-1}{t-2\ell}x^{t}.
\end{align*}
Therefore,
\begin{align*}
a(n,2m)&=\sum_{\ell=0}^{\lfloor \frac{n-3m}{2}\rfloor} \binom{m+\ell-1}{\ell}\binom{n-2m-2\ell-1}{m-1}.
\end{align*}
Similarly, we have
$$a(n,2m+1)=\sum_{\ell=0}^{\lfloor \frac{n-3m-1}{2}\rfloor} \binom{m+\ell-1}{\ell}\binom{n-2m-2\ell-1}{m}.$$
From these two identities we obtain the desired result. 
\end{proof}

Furthermore, we can use the Wilf-Zeilberger algorithm \cite{ABZ} to derive an additional recurrence relation for the sequence $a(n,m)$.

\begin{theorem}
For all $n, m\geq 0$,  
$$\left(m-n-2+\lfloor m/2\rfloor\right)a(n+2,m) +   \left(m-\lfloor m/2\rfloor\right)a(n+1,m) + na(n,m)=0$$
\end{theorem}
\begin{proof}
Let $F(n,\ell)$ be the expression
$$F(n,\ell):=\binom{m + \ell - 1}{\ell}\binom{n - m - \ell - 1}{n - m - \lfloor\frac{m}{2}\rfloor - \ell}(-1)^{n - m  -\lfloor\frac{m}{2}\rfloor - \ell}.$$
By the Wilf-Zeilberger's algorithm, we have that $F(n,\ell)$ satisfies this relation
\begin{multline*}
\left(m-n-2+\left\lfloor\frac m 2\right\rfloor\right)F(n+2,\ell) +   \left(m-\left\lfloor\frac m 2\right\rfloor\right)F(n+1,\ell) + nF(n,\ell)\\= G(n,\ell+1)-G(n,\ell),
\end{multline*}
with the certificate 
$$R(n,\ell)=\frac{\ell (-\ell - m + n)(-1 +\left\lfloor\frac m 2\right\rfloor)}{(-2 + \ell + m - n +\left\lfloor\frac m 2\right\rfloor)(-1 + \ell + m - n +\left\lfloor\frac m 2\right\rfloor)}.$$
That is, $R(n,\ell)=F(n,\ell)/G(n,\ell)$ is a rational function in both variables.     Summing over all $\ell$, the right-hand part cancels out, and we obtain the desired result.
\end{proof}

In the following corollary we give an asymptotic expression for the number of Arndt compositions with a fixed number of parts.  It is a direct application of transfer theorem (see Theorem 5.5 of \cite{Flajolet2}).

\begin{coro}
For a fixed positive integer $m$ we have 
$$a(n,m)\sim \frac{n^{m-1}}{2^{\lfloor m/2\rfloor}(m-1)!}.$$
\end{coro}
\begin{proof}
We know that 
\begin{align*}
a(n,2m)&=[x^n]\frac{x^{3m}}{(1-x)^{2m}(1+x)^m}=[x^n]\frac{f(x)}{(1-x)^{2m}},
\end{align*}
where $f(x)=x^{3m}/(1+x)^m$. From Theorem 5.5 of \cite{Flajolet2}, we conclude that 
\begin{align*}
a(n,2m)&\sim \frac{f(1)}{(2m-1)!}n^{2m-1}=\frac{n^{2m-1}}{2^m(2m-1)!}.
\end{align*}
If $m=2m+1$, then  $a(n,2m)\sim \frac{n^{2m}}{2^m(2m)!}$.
\end{proof}

Setting $y=1$ in Theorem \ref{gfun1} implies the following corollary. This is an alternative approach to this counting problem. Note that Hopkins and  Tangboonduangjit   \cite{HT1} proved this result by using a bijective approach. 
\begin{coro}\label{coro1a}
The generating function of the number of Arndt  compositions is
\begin{align*}
  A(x,1)&=\sum_{n\geq 0}a(n)x^n=\frac{1-x^2}{1-x-x^2}.
  \end{align*}
  Moreover,  $a(n)=F_n$ for all $n\geq1$.
\end{coro}

From Theorems \ref{teo2} and  \ref{teo2b} we obtain the following expressions (probably new) to calculate the Fibonacci number $F_n$ for $n\geq 1$:
\begin{align*}
F_n&=\sum_{m=0}^{n}\sum_{\ell=0}^{n-m-\lfloor\frac{m}{2}\rfloor}\binom{m + \ell - 1}{\ell}\binom{n - m - \ell - 1}{n - m - \lfloor\frac{m}{2}\rfloor - \ell}(-1)^{n - m  -\lfloor\frac{m}{2}\rfloor - \ell}\\
&=\sum_{m=0}^{n}\sum_{\ell=0}^{\lfloor (n-m-\lfloor\frac{m}{2}\rfloor)/2\rfloor}\binom{\lfloor\frac{m}{2}\rfloor + \ell - 1}{\ell}\binom{n - 2 \lfloor\frac{m}{2}\rfloor -2\ell-1}{\lfloor\frac{m-1}{2}\rfloor}.
\end{align*}

Let $p(n)$ denote the total number of parts (columns) over all compositions in $\Ar(n)$. The generating function for the sequence $p(n)$ is given by
\begin{align}\label{PPec}
P(x):=\sum_{n\geq 0}p(n)x^n&=\left.\frac{\partial A(x,y) }{\partial y}\right|_{y=1}=\frac{x(1 - x + x^3 - x^4)}{(1 - x - x^2)^2}\\&=x + x^2 + 3 x^3 + 6 x^4 + 11 x^5 + \bm{21} x^6 + 38 x^7+O(x^8).
\end{align}
For example, from Figure \ref{Fig1} we can verify that the number of parts (columns) of all compositions  in $\Ar(6)$ is 21.

We need the following  for the asymptotics of linear recurrences relations (see \cite{Flajolet2}). Assume that a rational generating function $f(x)/g(x)$, with $f(x)$ and $g(x)$ relatively prime and $g(0)\neq 0$, has a unique pole $1/\beta$ with the smallest modulus. Then, if the multiplicity of $1/\beta$ is $\nu$, we have
\begin{align}\label{FlaTheo}
[x^n]\frac{f(x)}{g(x)}\sim \nu \frac{(-\beta)^\nu f(1/\beta)}{g^{(\nu)}(1/\beta)}\beta^n n^{\nu-1}.
\end{align}
From \eqref{PPec} and  \eqref{FlaTheo}  we obtain that 
\begin{align}\label{aspn} 
p(n)\sim \frac{3-\sqrt 5}{5}\left(\frac{1+\sqrt 5}{2} \right)^nn.
\end{align}
\begin{theorem}
The expected number of parts in $\Ar(n)$ is asymptotically
$$\left(\frac{3}{\sqrt 5}-1\right)n.$$
\end{theorem}
\begin{proof}
 Let $X_{\Arr}(n)$ denote the number of parts in a  random Arndt composition in $\Ar(n)$.
The expected value is given by
\begin{align}
\E[X_{\Arr}(n)]=\frac{1}{a(n)}\sum_{\sigma \in\Ar(n)}  \parts(\sigma)=\frac{p(n)}{F_n}.
\end{align}
From  $F_n\sim \frac{1}{\sqrt 5}\left(\frac{1+\sqrt 5}{2} \right)^n$ and \eqref{aspn} we conclude the desired result.
\end{proof}

\subsection{Anti-palindromic compositions}

A  composition $\sigma=(\sigma_1, \sigma_2,\ldots, \sigma_\ell)$  of a positive integer $n$ is \emph{anti-palindromic} if $\sigma_i\neq \sigma_\ell+1-i$ for all $i\neq (\ell+1)/2$.  This family of compositions was recently studied by Andrews, Just, and Simay \cite{AJS}. For example, $(\texttt{1}, \texttt{2}, \texttt{6}, \texttt{3}, \texttt{2})$ is an anti-palindromic composition of 14. Let $\mathcal{AP}$ denote the set of all anti-palindromic compositions. 
For each anti-palindromic composition of $n$ with  $\ell$ parts, it is possible to form $2^{\lfloor \ell/2\rfloor}$ \emph{flip-equivalent} anti-palindromic compositions of $n$ with $\ell$ parts by switching any number of the pairs $\sigma_i$ and $\sigma_{\ell-i+1}$ ($i\neq (\ell+1)/2$). The sets of flip-equivalent anti-palindromic compositions of $n$ form a partition of the set $\mathcal{AP}$ and each equivalence class is called a \emph{reduced anti-palindromic composition} of $n$ with $\ell$ parts. For example, the following anti-palindromic compositions are flip-equivalent
$$(\texttt{1}, \texttt{2}, \texttt{6}, \texttt{3}, \texttt{2}), \quad (\texttt{2}, \texttt{2}, \texttt{6}, \texttt{3}, \texttt{1}), \quad (\texttt{1}, \texttt{3}, \texttt{6}, \texttt{2}, \texttt{2}), \quad \text{and} \quad (\texttt{2}, \texttt{3}, \texttt{6}, \texttt{2}, \texttt{1}).$$
Let $\mathcal{RAP}$ be the set of all reduced anti-palindromic compositions. We introduce the following bivariate generating functions with respect to the number of parts and weight over the sets $\mathcal{AP}$ and $\mathcal{RAP}$, that is
$$Ap(x,y):=\sum_{\sigma\in \mathcal{AP}}x^{|\sigma|}y^{\parts(\sigma)} \quad \text{and} \quad Bp(x,y):=\sum_{\sigma\in \mathcal{RAP}}x^{|\sigma|}y^{\parts(\sigma)}.$$

\begin{theorem}\label{gfun1ap}
The generating function for anti-palindromic compositions  with respect to the number of
parts and weight is given by
$$Ap(x,y)=\frac{1 - x - x^2 + x^3 + x y - x^3 y}{1 - x - x^2 + x^3 - 2 x^3 y^2}.$$
\end{theorem}

\begin{proof}
Let $\sigma=(\sigma_1, \sigma_2, \dots, \sigma_{2m})$ be an anti-palindromic composition of $[n]$ with $2m$ parts. From the definition we have the condition $\sigma_i\neq \sigma_{2m+1-i}$, for all $i\neq (2m+1)/2$. Notice that the columns $i$-th and $(2m+1-i)$-th contribute to the generating function the term 
\begin{align}\label{ecflip}
y^2\sum_{i\geq 1} x^i \sum_{j\geq 1, j\neq i} x^j=y^2\sum_{i\geq 1} x^i\left(\frac{x}{1-x}-x^i\right)=\frac{2x^3}{(1-x)^2 (1 + x)}.
\end{align}
 Therefore the composition $\sigma$ contributes to the generating function the term $$\sum_{m\geq 0}\left(\frac{2x^3}{(1-x)^2 (1 + x)}\right)^my^{2m}=\frac{(1-x)^2(1+x)}{1 - x - x^2 + x^3 - 2 x^3 y^2}.$$

If the number of parts is odd, then from a similar argument the contribution to the generating function is given by $$\frac{x(1-x^2)y}{1 - x - x^2 + x^3 - 2 x^3 y^2}.$$
Summing the above two generating functions, we obtain the desired result.
\end{proof}

Notice that if we divide \eqref{ecflip} by 2, then we obtain the generating function for the reduced anti-palindromic compositions. 
\begin{coro}
The generating function for reduced anti-palindromic compositions with respect to the number of
parts and weight is given by
$$Bp(x,y)=\frac{1 - x - x^2 + x^3 + x y - x^3 y}{1 - x - x^2 + x^3 - x^3 y^2}.$$
\end{coro}

Since $Bp(x,y)=A(x,y)$ (see Theorem \ref{gfun1}) we obtain the following result.
\begin{theorem}\label{bijection}
The number of Arndt compositions of $n$ with $m$ parts is equal to the number of reduced anti-palindromic compositions $n$ with $m$ parts.
\end{theorem}

It is also possible to establish an explicit bijection between these sets.  First, for convenience, we will take as a representative of a reduced anti-palindromic composition one that satisfies the condition $\sigma_i>\sigma_{\ell-i+1}$ for each pair of summands with $i<\frac{\ell+1}{2}$. For example, for the compositions given above, the representative is $(\texttt{2}, \texttt{3}, \texttt{6}, \texttt{2}, \texttt{1})$. Once we impose this condition, the connection between Arndt compositions and reduced anti-palindromic compositions is immediate, as in both cases, we are comparing pairs of summands, the first greater than the other. In the case of Arndt compositions, the pairs of summands are ordered consecutively and in the case of reduced anti-palindromic compositions they are ordered at each side of the composition.

For example, the representation of $(\texttt{2}, \texttt{3}, \texttt{6}, \texttt{2}, \texttt{1})$ as an Arndt composition would be $(\texttt{2}, \texttt{1}, \texttt{3}, \texttt{2}, \texttt{6})$. To go from one to the other, we take each pair of summands at the sides and reorder them to be adjacent in an Arndt composition. When the length of the reduced anti-palindromic composition is odd, the summand in the middle (the one that does not need to be compared to another) becomes the last in the Arndt composition. It is clear that this algorithm defines a bijection between both families of compositions.

\section{The Last Part}

We denote by $\last(\sigma)$ the weight of the last part in a composition $\sigma$.  We introduce the bivariate generating function to count the number of Arndt compositions with respect to the
weight and the last part:
$$B(x,y):=\sum_{\sigma\in \Ar}x^{|\sigma|}y^{\last(\sigma)}.$$
The coefficient $x^ny^m$ of $B(x,y)$ is equal to the number of Arndt compositions of $n$ whose last part is equal to $m$. In Theorem \ref{gfun4} we give a rational generating function for  $B(x,y)$.

\begin{theorem}\label{gfun4}
The generating function for Arndt compositions  with respect to the weight and the  last part  is given by
$$B(x,y)=\frac{1 - x - x^2 - x^2 y + 2 x^3 y + 2 x^4 y - x^5 y - x^4 y^2}{ (1 - x - x^2) (1 - x y) (1 - x^2 y)}.$$
\end{theorem}
\begin{proof}
Let $\sigma=(\sigma_1, \sigma_2)$  be an Arndt composition with two parts.  If $\sigma_1=j\geq 2$, then  the bivariate generating function is  
$$\sum_{j\geq 2}x^{j}\left(xy+(xy)^2+\cdots +(xy)^{j-1}\right)=\sum_{j\geq 2}x^j\frac{xy-(xy)^j}{1-xy}=\frac{x^3y^2}{(1-x)(1+x^2y)}.$$
Therefore the generating function for Arndt compositions with an even number of parts with respect to the last part is given by
 $$1+\sum_{m\geq 1}\left(\frac{x^3}{(1-x)(1+x^2)}\right)^{m-1}\frac{x^3y^2}{(1-x)(1+x^2y)}=1+\frac{x^3 (1 - x^2) y}{(1 - x - x^2) (1 - x^2 y)}.$$
 Analogously, if the number of parts is odd, then generating functions is
 $$\sum_{m\geq 0}\left(\frac{x^3}{(1-x)(1+x^2)}\right)^m\frac{xy}{1-xy}=\frac{(1 - x)^2 x (1 + x) y}{(1 - x - x^2) (1 - x y)}.$$
Adding the last two equations we obtain the desired result.
\end{proof}

\begin{figure}[h]
\centering
  \includegraphics[scale=0.8]{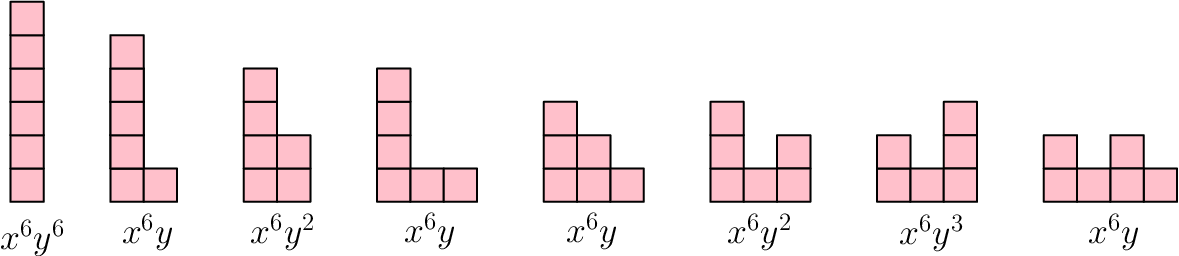}
  \caption{Weights for Arndt compositions (bar graphs) in $\Ar(3)$.}
  \label{Fig2}
\end{figure}

As a series expansion, the generating function $B(x,y)$ begins with
\begin{multline*}
B(x,y)=1+x y+x^2 y^2+x^3 \left(y^3+y\right)+x^4 \left(y^4+2
   y\right)\\+x^5 \left(y^5+2 y^2+2 y\right)+\bm{x^6
   \left(y^6+y^3+2 y^2+4 y\right)}+O\left(x^7\right).
   \end{multline*}
Figure \ref{Fig2} shows the weights of the  Arndt compositions   corresponding to the bold coefficient in the above series.

Let $b(n,m)$ denote the number of Arndt compositions of $n$, whose last part is equal to $m$. Note that $a(n)=\sum_{m\geq 1}b(n,m)$. In Table \ref{tab2} we show the first few values of the sequence $b(n,m)$. 
\begin{table}[htp]
\begin{center}
\begin{tabular}{c|cccccccccccccccc}
$n\backslash m$ &0& 1 & 2 & 3 & 4 & 5 & 6 & 7 & 8 & 9 & 10     \\ \hline
 0 &1 &  &  &  &  & &  &  &  &  &  \\
1 & 0 & 1 &  &  &  &  &  &  &  &  &  \\
2 & 0 & 0 & 1 &  &  &  &  &  &  &  &  \\
3 & 0 & 1 & 0 & 1 &  &  &  &  &  &  &  \\
4 & 0 & 2 & 0 & 0 & 1 &  & &  &  &  &  \\
5 & 0 & 2 & 2 & 0 & 0 & 1 &  &  &  &  &  \\
6 & 0 & 4 & 2 & 1 & 0 & 0 & 1 &  &  &  &  \\
7 & 0 & 6 & 3 & 2 & 1 & 0 & 0 & 1 &  &  &  \\
8 & 0 & 10 & 5 & 3 & 1 & 1 & 0 & 0 & 1 &  &  \\
9 & 0 & 16 & 8 & 4 & 3 & 1 & 1 & 0 & 0 & 1 &  \\
10 & 0 & 26 & 13 & 7 & 4 & 2 & 1 & 1 & 0 & 0 & 1 \\
\end{tabular}
\end{center}
\caption{Values of $b(n,m)$, for $0\leq n \leq 10$ and $0\leq m \leq10$.}
\label{tab2}
\end{table}
\begin{theorem}\label{teo2bc}
For all $n\geq  1$ we have 
$$b(n,m)=\begin{cases}
t_1(n-m), & \text{if} \ m\leq n <2m;\\ 
t_1(n-m) + t_2(n-2m), & \text{if} \ 2m\leq n;\\ 
0, & \text{otherwise;}
\end{cases}$$
where $t_1(0)=1$, $t_1(1)=0$, $t_1(n)=F_{n-2}$  for $n\geq 2$ and  $t_2(0)=0$, $t_2(1)=1$,$t_2(n)=F_{n-1}$ for $n\geq 2$.
\end{theorem}
\begin{proof}
The generating function $B(x,y)$ can be rewritten as
$$B(x,y)=-F(x) + T_1(x,y) + T_2(x,y),$$
where $F(x)$ is the generating function of the Fibonacci numbers, that is $F(x)=x/(1-x-x^2)$. Moreover, 
$$T_1(x,y)=\frac{(1 - x)^2 (1 + x)}{(1 - x - x^2) (1 - x y)} \quad \text{and} \quad T_2(x,y)=\frac{(1 - x) x (1 + x)}{(1 - x - x^2) (1 - x^2 y)}.$$
By standard  methods the  generating functions of the sequences $t_1(n)$ and $t_2(n)$, defined in the statement of this theorem, are given by
$$T_1(x)=\sum_{n\geq 0}t_1(n)x^n=\frac{(1 - x)^2 (1 + x)}{(1 - x - x^2)} \quad \text{and} \quad T_2(x)=\sum_{n\geq 0}t_2(n)x^n=\frac{(1 - x)x(1 + x)}{(1 - x - x^2)}.$$ 
By the Cauchy product we have
\begin{align*}
T_1(x,y)&=T_1(x) \frac{1}{1-xy}=\sum_{n\geq 0}\sum_{m=0}^nt_1(n-m)y^mx^n=\sum_{m\geq 0}y^m\sum_{n\geq 0}t_1(n)x^{n+m}\\
T_2(x,y)&=T_2(x) \frac{1}{1-x^2y}=\sum_{n\geq 0}\sum_{m=0}^nt_2(n-m)y^mx^{2n}=\sum_{m\geq 0}y^m\sum_{n\geq 0}t_2(n)x^{n+2m}.
\end{align*}
By comparing the coefficients we obtain the desired result. 
\end{proof}

Notice that for all  $n\geq 2m+2$, we have the equality  
$$b(n,m)=F_{n-m-2} + F_{n-2m-1}.$$
For example, the number of Arndt compositions whose last part is equal to \texttt{1} is given by $b(n,1)=2F_{n-3}$ for all $n\geq 4$. Our result  generalizes the Corollary 2.2 of \cite{HT1}.

Let $b^{(\leq m)}(n)$ and  $b^{(\geq m)}(n)$  be the number of Arndt composition whose last part is at most $m$ and at least $m$, respectively. 

\begin{coro}\label{cor:lastsumatm} 
For $k\geq 1$ and  $n\geq 2k+2$, we have  
$$b^{(\leq k)}(n)=F_{n}- F_{n-k-1} - F_{n-2k-2} \quad \text{and} \quad b^{(\geq k)}(n)=F_{n-k}- F_{n-2k}.$$ 
\end{coro}
\begin{proof}
    From the definitions we have
        \begin{align*}
        b^{(\leq k)}(n)&=\sum_{j=1}^kb(n,j)=\sum_{j=1}^k(F_{n-j-2} + F_{n-2j-1})\\
        &=\sum_{j=1}^k((F_{n-j}-F_{n-j-1}) + (F_{n-2j}-F_{n-2j-2}))\\
        &=F_{n-1} - F_{n-k-1} + F_{n-2} - F_{n-2k-2}\\
        &=F_n- F_{n-k-1} - F_{n-2k-2}. 
    \end{align*}
 The second identity can be obtained in a similar way.
 \end{proof}

From Corollary \ref{cor:lastsumatm} there are $F_{n-2}+F_{n-4}$ Arndt compositions whose last part is greater than 1. This result  generalizes the Corollary 2.2 of \cite{HT1}.

 Finally, we also have an asymptotic formula for $b(n,m)$.

 \begin{theorem}
For a fixed positive integer $m$ we have 
$$b(n,m)\sim\frac{\varphi^{n-m-2}}{\sqrt 5}\left(1+ \frac{1}{\varphi^{m-1}}\right),$$
    where $\varphi=\frac{1+\sqrt{5}}{2}$ is the golden ratio.
   \end{theorem}
\begin{proof}
From Theorem \ref{teo2bc} we know that for $2m\leq n$
$$b(n,m)=t_1(n-m) + t_2(n-2m)= F_{n-m-2} + F_{n-2m-1}.$$
From  $F_n\sim \frac{1}{\sqrt 5}\varphi^n$, we have
$$b(n,m) \sim   \frac{1}{\sqrt 5}\left(\varphi^{n-m-2}+ \varphi^{n-2m-1}\right)=\frac{\varphi^{n-m-2}}{\sqrt 5}\left(1+ \varphi^{-m+1}\right).$$ 
\end{proof}

Let $d(n)$ denote the sum of the weight of the last part over all compositions in $\Ar(n)$. The generating function for the sequence $d(n)$ is given by
\begin{align*}
D(x):=\sum_{n\geq 0}d(n)x^n&=\left.\frac{\partial B(x,y) }{\partial y}\right|_{y=1}=\frac{x(1 + x - x^3)}{1 - x - 2 x^2 + x^3 + x^4}\\&=x + 2 x^2 + 4 x^3 + 6 x^4 + 11 x^5 + \bm{17}x^6 + 29 x^7 +O(x^8).
\end{align*}
From Figure \ref{Fig1} we can verify that the sum of the last parts (last column) of all compositions in $\Ar(6)$ is 17. Notice that the sequence $d(n)$ corresponds to the sequence \seqnum{A014217}, therefore
$$d(n)=\left\lfloor\left(\frac{1+\sqrt 5}{2}\right)^n\right\rfloor, \quad n\geq 1.$$
Our results leave  a new combinatorial interpretation for this sequence.

 \begin{coro}
The expected number of the weight of the last part in $\Ar(n)$ is asymptotically $\sqrt{5}\approx 2.23607$.
\end{coro}

\section{First generalization of Arndt compositions}

Hopkins and Tangboonduangjit \cite{HT1} generalized Arndt compositions by requiring a greater decrease or increase. Specifically, given a positive integer $k$, let $\Ar_k(n)$ denote the set of Arndt compositions of $n$ such that $\sigma_{2i-1}>\sigma_{2i}+k$ for each positive integer $i$. The elements in $\Ar_k(n)$ will be called \emph{$k$-Arndt compositions}. Furthermore, we will denote by
$a_k(n,m)$ the number of such compositions of $n$ that have exactly $m$ parts and $a_k(n):=\sum_{m\geq 1}a_k(n,m)$. For example, $a_3(10)= 10$, the relevant compositions being
\begin{align*}
&(\texttt{10}), \quad (\texttt{9}, \texttt{1}), \quad (\texttt{8}, \texttt{2}), \quad (\texttt{8}, \texttt{1}, \texttt{1}),  \quad (\texttt{7}, \texttt{3}), \quad (\texttt{7}, \texttt{2}, \texttt{1}),\\
&(\texttt{7}, \texttt{1}, \texttt{2}), \quad (\texttt{6}, \texttt{1}, \texttt{3}), \quad (\texttt{6}, \texttt{2}, \texttt{2}), \quad (\texttt{5}, \texttt{1}, \texttt{4}).
\end{align*}
   
   Let $\Ar_k$ denote the set of all $k$-Arndt  compositions.   We introduce a bivariate generating function to count the number of $k$-Arndt compositions with respect to the weight and number of parts:
$$A_k(x,y):=\sum_{\sigma\in \Ar_k}x^{|\sigma|}y^{\parts(\sigma)}.$$

We have the following generating function formulas for the preceding sequences.
\begin{theorem}\label{gfung1}
The generating function for the $k$-Arndt compositions  with respect to the number of
parts and weight is given by
\begin{align*}A_k(x,y)=\begin{cases}
\dfrac{(1 - x^2) (1 - x (1 - y))}{1 - x - x^2 + x^3 - x^{3 + k} y^2},& \text{if } k\geq 0;\\
\dfrac{(1 - x^2) (1 - x (1 - y))}{1 - x - x^2 (1 + y^2) + x^3 (1 - y^2) + y^2 x^{2 - k}},& \text{if } k< 0;
\end{cases}
\end{align*}
\end{theorem}
\begin{proof}
If $k\geq 0$, then the  $k$-Arndt composition with two parts $(\sigma_1,\sigma_2)$ ($j=\sigma_1>\sigma_2+k$) are enumerated by
$$\sum_{j\geq 2 + k}x^{j}y\left(x+x^2+\cdots +x^{j-1-k}\right)y=\frac{x^{3+k}y^2}{(1-x)^2(1+x)}.$$
Therefore the generating function for the $k$-Arndt compositions with an even number of parts is
 $$\sum_{m\geq 0}\left(\frac{x^{3+k}}{(1-x)^2(1+x)}\right)^my^{2m}=\frac{(1 - x)^2 (1 + x)}{1 - x - x^2 + x^3 - x^{3 + k} y^2}.$$
 Analogously, if the number of parts is odd, then the generating function is
  $$\sum_{m\geq 0}\left(\frac{x^{3+k}}{(1-x)^2(1+x)}\right)^my^{2m}\frac{xy}{1-x}=\frac{(1 - x^2)xy}{1 - x - x^2 + x^3 - x^{3 + k} y^2}.$$
Adding the last two equations we obtain the desired result. 
The proof for case $k<0$ follows from a similar argument.
\end{proof}

As a series expansion, the generating functions $A_3(x,y)$ and $A_{-3}(x,y)$begins with
\begin{align*}
A_3(x,y)&=1 + y x + y x^2 + y x^3 + y x^4 +  y x^5 + (y + y^2) x^6 + (y + y^2 + y^3) x^7+O\left(x^8\right),\\
A_{-3}(x,y)&=1 + y x + (y + y^2) x^2 + (y + 2 y^2 + y^3) x^3 + (y + 3 y^2 + 3 y^3 +
     y^4) x^4 \\
     & \ + (y + 3 y^2 + 6 y^3 + 4 y^4 + y^5) x^5 + (y + 4 y^2 + 
    9 y^3 + 10 y^4 + 5 y^5 + y^6) x^6 +O\left(x^7\right).
   \end{align*}

Setting $y=1$ in Theorem \ref{gfung1} implies the following corollary.
\begin{coro}\label{coro1ga}
The generating function of the number of $k$-Arndt  compositions is
\begin{align*}A_k(x)=\begin{cases}
\dfrac{1 - x^2}{1 - x - x^2 + x^3 - x^{k+3}},& \text{if } k\geq 0;\\
\dfrac{1 - x^2}{1 - x - 2 x^2 + x^{2 - k}},& \text{if } k< 0;
\end{cases}
\end{align*}
\end{coro}
Note that if $k=0$, we recover the Corollary \ref{coro1a}.

\section{A second generalization of Arndt compositions}

The goal of this section is to introduce a second generalization of Arndt compositions. Specifically, given a positive integer $k$, let $\Ar^{(k)}(n)$ denote the set of compositions of $n$ such that $\sigma_{ki-(k-1)}>\sigma_{ki-(k-2)}>\cdots > \sigma_{ki}$ for each positive integer $i$. The elements in $\Ar^{(k)}(n)$ will be called \emph{$k$-block Arndt compositions}. Furthermore, we will denote by
$a^{(k)}(n,m)$ the number of such compositions of $n$ that have exactly $m$ parts and $a^{(k)}(n):=\sum_{m\geq 1}a^{(k)}(n,m)$. For example, $a^{(3)}(10)= 18$, the relevant compositions being
\begin{align*}
&(\texttt{10}), \quad (\texttt{9}, \texttt{1}), \quad (\texttt{8}, \texttt{2}), \quad (\texttt{7},\texttt{3}), \quad (\texttt{7},\texttt{2},\texttt{1}), \quad (\texttt{6},\texttt{4}), \quad (\texttt{6},\texttt{3},\texttt{1}), \quad (\texttt{6},\texttt{2},\texttt{1},\texttt{1}),\\
&(\texttt{5},\texttt{4},\texttt{1}), \quad (\texttt{5},\texttt{3},\texttt{2}), \quad (\texttt{5},\texttt{3},\texttt{1},\texttt{1}), \quad (\texttt{5},\texttt{2},\texttt{1},\texttt{2}), \quad (\texttt{4},\texttt{3},\texttt{2},\texttt{1}), \quad (\texttt{4},\texttt{3},\texttt{1},\texttt{2}), \\
&(\texttt{4},\texttt{2},\texttt{1},\texttt{3}), \quad (\texttt{3},\texttt{2},\texttt{1},\texttt{4}), \quad (\texttt{4},\texttt{2},\texttt{1},\texttt{2},\texttt{1}), \quad (\texttt{3},\texttt{2},\texttt{1},\texttt{3},\texttt{1}).
\end{align*}

Let $\Ar^{(k)}$ denote the set of all $k$-block Arndt compositions.  Consider the bivariate generating function
$$A^{(k)}(x,y):=\sum_{\sigma\in \Ar^{(k)}}x^{|\sigma|}y^{\parts(\sigma)}.$$

\begin{theorem}
The generating function for the $k$-block Arndt compositions with respect to the number of
parts and weight is given by
\begin{align*}
A^{(k)}(x,y)=\frac{1}{1-J_k(x,y)}\left(\sum_{j=0}^{k-1}J_j(x,y)\right),\end{align*}
where $J_j(x,y)=x^{\binom{j+1}{2}}y^j\prod_{\ell=1}^j\frac{1}{1-x^{\ell}}$ for $j\geq 1$ and $J_0(x,y)=1$.
\end{theorem}
\begin{proof}
Let $\beta=(\beta_1, \beta_2, \dots, \beta_j)$  be a $k$-block Arndt composition with $j$ parts $(j\leq k)$. This kind of compositions are integer partitions with exactly $j$ distinct parts. The generating functions for this family of partitions is given by (cf. \cite[pp. 105]{Comtet})
$$J_j(x,y)=y^j\prod_{\ell=1}^j\frac{x^\ell}{1-x^{\ell}}=x^{\binom{j+1}{2}}y^j\prod_{\ell=1}^j\frac{1}{1-x^{\ell}}.$$
The  $k$-block Arndt compositions are the concatenation of blocks of $k$ columns in decreasing order and at most $j\leq k-1$ additional columns in decreasing order; see Figure \ref{gencase} for a pictorial description.
 \begin{figure}[H]
\centering
  \includegraphics[scale=0.7]{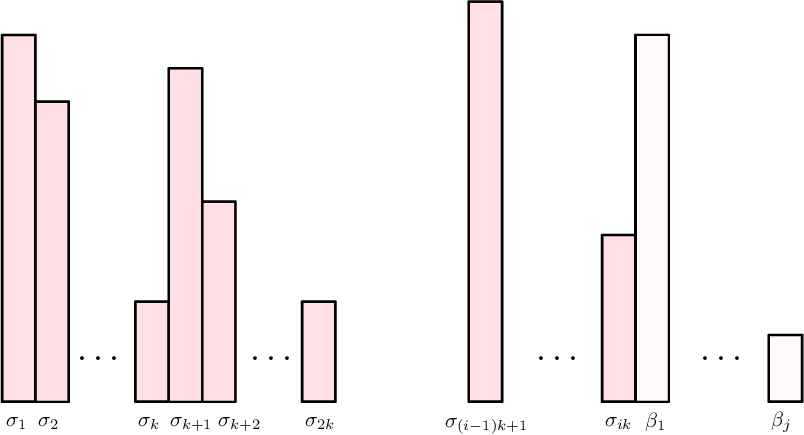}
  \caption{Decomposition of a $k$-block Arndt composition.}
  \label{gencase}
\end{figure}
Therefore, the generating function is given by
 $$ \sum_{i\geq 0}(J_k(x,y))^i\left(\sum_{j=0}^{k-1}J_j(x,y)\right)=\frac{1}{1-J_k(x,y)}\left(\sum_{j=0}^{k-1}J_j(x,y)\right).
 $$
 \end{proof}
 For example, for $k=3$ and $k=4$ we obtain the bivariate generating functions:
 \begin{align*}
 A^{(3)}(x,y)&=\frac{(1-x^3) (1 - x - x^2 + x^3 + x y - x^3 y + x^3 y^2)}{1 - x - x^2 + x^4 + x^5 - x^6(1 + y^3)},\\
  A^{(4)}(x,y)&=\frac{(1 - x^4)(1 - x - x^2 + x^4 + x^5 - x^6 + x y - 
   x^3 y - x^4 y + x^6 y + x^3 y^2 - x^6 y^2 + x^6 y^3)}{1 - x - x^2 + 2 x^5 - x^8 - x^9 + x^{10} - x^{10}y^4},
 \end{align*}

Setting $y=1$ in the above generating functions we obtain, respectively,
the generating functions of the number of $3$-block and $4$-block Arndt compositions:

 \begin{align*}
 A^{(3)}(x,1)&=\frac{1 - x^2 + x^5 - x^6}{1 - x - x^2 + x^4 + x^5 - 2x^6}\\
 &=1 + x + x^2 + 2 x^3 + 2 x^4 + 3 x^5 + 4 x^6 + 6 x^7 + 8 x^8 + 
 13 x^9 + O(x^{10}),\\
  A^{(4)}(x,1)&=\frac{(1 - x) (1 + x) (1 + x^2) (1 - x^2 + x^5)}{1 - x - x^2 + 2 x^5 - x^8 - x^9}\\
  &=1 + x + x^2 + 2 x^3 + 2 x^4 + 3 x^5 + 4 x^6 + 5 x^7 + 6 x^8 + 8 x^9 + 
 10 x^{10} + O(x^{11}).
 \end{align*}

For example, these are some of the $3$-block Arndt compositions of $14$ with $6$ parts
$$(\texttt{5}, \texttt{2}, \texttt{1}, \texttt{3}, \texttt{2}, \texttt{1}), 	\quad (\texttt{3}, \texttt{2}, \texttt{1}, \texttt{5}, \texttt{2}, \texttt{1}), \quad (\texttt{4}, \texttt{2}, \texttt{1}, \texttt{4}, \texttt{2}, \texttt{1}), \quad (\texttt{4}, \texttt{3}, \texttt{1},
  \texttt{3}, \texttt{2}, \texttt{1}), \quad (\texttt{3}, \texttt{2}, \texttt{1}, \texttt{4}, \texttt{3}, \texttt{1}).$$

\section{Acknowledgments}
The authors were partially supported by Universidad Nacional de Colombia, Project No. 57340.


\begin{thebibliography}{20}

\bibitem{AJS} G. E.~Andrews, M.~Just, and G.~Simay, Anti-palindromic compositions, Fibonacci Quart.
\textbf{60} (2022), 164--176.

\bibitem{Comtet}
L.~Comtet. Advanced Combinatorics.The Art of Finite and Infinite Expansions,  D. Reidel Publishing Company, Dordrecht, Holland 1974.

\bibitem{CM} R.~Graham, D.~Knuth, and O.~Patashnik, Concrete Mathematics: A Foundation for Computer Science.   Addison-Wesley Professional; 2nd edition, 1994.

\bibitem{HT1} B.~Hopkins and A.~Tangboonduangjit, Verifying and generalizing Arndt’s compositions, Fibonacci
Quart. \textbf{60}(5) (2022), 181--186.

\bibitem{HT2} B.~Hopkins and A.~Tangboonduangjit, Arndt and De Morgan integer compositions. arXiv:2307.12434v1 (2023).

\bibitem{ABZ} M. Petko\v{v}sek, H. Wilf, and D. Zeilberger, \emph{A=B}, \emph{A. K. Peters}, Ltd. 1996.

\bibitem{Flajolet2} 
R.~Sedgewick and P.~Flajolet, An Introduction to the Analysis of Algorithms, 2nd Edition, Addison-Wesley, 2013.

\bibitem{OEIS} N.~J.~A.~Sloane. The On-Line Encyclopedia of Integer Sequences, \url{http://oeis.org/}.

\end{thebibliography}
\end{document}